\documentclass[11pt,a4paper,leqno,twoside, headinclude]{amsart}

\title[Non-formal Homogeneous Spaces]{Non-formal Homogeneous Spaces}
\def\titl{Non-formal Homogeneous Spaces}
\def\auth{Manuel Amann}
\date{May 18th, 2012}
\keywords{\noindent  homogeneous space, non-formal manifold}
\thanks{The author was supported by a Research Grant of the German Research Foundation.}
\subjclass[2010]{57N65 (Primary), 57T15, 57T20 (Secondary)}
\author{\auth}

\usepackage[english]{babel}

\usepackage[colorlinks,pdftex, plainpages=false]{hyperref}
\hypersetup{pdftitle=\titl, pdfauthor=\auth, pdftoolbar=false,
plainpages=false, hyperindex=true, pdfdisplaydoctitle=true}

\hypersetup{colorlinks,%
citecolor=black,%
filecolor=black,%
linkcolor=black,%
urlcolor=black,%
pdftex}

\usepackage{color}

\usepackage{amsmath, amssymb, amscd, amsthm}
\usepackage{stmaryrd}
\usepackage{latexsym}
\usepackage[all]{xy}
\usepackage{pb-diagram}
\usepackage{rotating}
\usepackage{multicol}
\usepackage{lscape}
\usepackage{wasysym}
\usepackage{longtable}

\xyoption{all}

\newtheorem{theo}{Theorem}[section]
\newtheorem{main}{Theorem}
\newtheorem*{main*}{Theorem}
\newtheorem{mainprop}[main]{Proposition}
\newtheorem*{mainprop*}{Proposition}

\newtheorem{prop}[theo]{Proposition}
\newtheorem{defi2}[theo]{Definition}
\newtheorem*{defi2*}{Definition}

\newenvironment{defi*}{\begin{defi2*}\normalfont}{\end{defi2*}}

\newenvironment{defin*}[1]{\begin{defi2*}[#1]\normalfont}{\end{defi2*}}
\newtheorem{rem2}[theo]{Remark}
\newenvironment{rem}{\begin{rem2}\normalfont}{\hfill$\boxbox$\end{rem2}}
\newtheorem{lemma}[theo]{Lemma}
\newtheorem{cor}[theo]{Corollary}

\newtheorem{ques}[theo]{Question}
\newtheorem*{conj*}{Conjecture}
\newtheorem*{theo*}{Theorem}
\newtheorem*{ques*}{Question}
\newtheorem*{mi2}{Main Idea}

\newtheorem{ex2}[theo]{Example}
\newenvironment{ex}{\begin{ex2}\normalfont}{\hfill$\boxbox$\end{ex2}}
\newtheorem{exer2}[theo]{Exercise}

\newtheorem{alg2}[theo]{Algorithm}

\newcommand{\cc}{{\mathbb{C}}}                                     
\newcommand{\hh}{{\mathbb{H}}}                                     
\newcommand{\nn}{{\mathbb{N}}}                                     
\newcommand{\qq}{{\mathbb{Q}}}                                     
\newcommand{\rr}{{\mathbb{R}}}                                     
\newcommand{\s}{{\mathbb{S}}}                                      
\newcommand{\zz}{{\mathbb{Z}}}                                     
\newcommand{\SO}{{\mathbf{SO}}}                                    
\newcommand{\U}{{\mathbf{U}}}                                      
\newcommand{\SU}{{\mathbf{SU}}}                                    
\newcommand{\Sp}{{\mathbf{Sp}}}                                    
\newcommand{\B}{{\mathbf{B}}}                                      
\newcommand{\E}{{\mathbf{E}}}                                      
\newcommand{\dif} {{\operatorname{d}}}                             
\newcommand{\In} {{\,\subseteq\,}}                                 
\newcommand{\im} {{\operatorname{im\,}}}                           
\newcommand{\id}{{\operatorname{id}}}                              
\newcommand{\APL}{{\operatorname{A_{PL}}}}                         
\newcommand{\rk}{{\operatorname{rk\,}}}                            
\newcommand{\W}{{\operatorname{W}}}                                
\newcommand{\co}{\colon\thinspace}                                 
\newcommand{\even}{\textrm{even}}                                  
\newcommand{\odd}{\textrm{odd}}                                    
\newcommand{\vproof}{{\begin{flushright} \qed                      
                      \end{flushright}}}

\newcommand{\comment}[1]{}                                         
\newcommand{\xto}[1]{\xrightarrow{#1}}                             
\newcommand{\hto}[1]{\overset{#1}{\hookrightarrow}}                
\newcommand{\biq}[2]{#1\;\!\!\!\sslash \;\!\!\!#2}                 
\newcommand{\ack}{\noindent\textbf{Acknowledgements. }}            
\newcommand{\str}{\noindent\textbf{Structure of the article. }}    

\newenvironment{prf}{\begin{proof}[\textsc{Proof}]} {\end{proof}}     

\begin{document}

\maketitle \thispagestyle{empty}


\begin{abstract}
Several large classes of homogeneous spaces are known to be formal---in the sense of Rational Homotopy Theory. However, it seems that far fewer examples of non-formal homogeneous spaces are known.

In this article we provide several construction principles and characterisations for non-formal homogeneous spaces, which will yield a lot of examples. This will enable us to prove that, from dimension $72$ on, such a space can be found in each dimension.
\end{abstract}


\section*{Introduction}
Homogeneous spaces form a very well-studied and interesting class of manifolds. They appear abundantly in geometry and topology.
Our focus will lie on their topological properties adopting the viewpoint of Rational Homotopy Theory.
Also from this perspective homogeneous spaces bear remarkable properties, which they share with the larger class of biquotients they are contained in.

The group $G$ will always be taken to be a compact connected Lie group and let
$H\In G\times G$ be a closed Lie subgroup. Then $H$ acts on $G$ on
the left by $(h_1,h_2)\cdot g=h_1gh_2^{-1}$. The orbit space of this
action is called the \emph{biquotient} $\biq{G}{H}$ of $G$ by $H$. If the action of $H$ on $G$ is free, then
$\biq{G}{H}$ possesses a manifold structure. This is the only case
we shall consider. Clearly, the category of biquotients contains the
one of homogeneous spaces; in this special case the inclusion of $H$ into the first $G$ factor is trivial.

It is a classical result that biquotients are \emph{rationally elliptic spaces}. Moreover, they admit what is called a \emph{pure model}---see \cite{FHT01}, p.~435, for the definition and examples \cite{FHT01}.32.2, p.~448 and \cite{FHT01}.15.1, p.~218 for the homogeneous case, respectively \cite{Kap}.1, p.~2 for general biquotients. This makes them an interesting yet manageable source of examples within the realm of Rational Homotopy Theory and they definitely do constitute a field of study well-worth
the attention it is paid.

The topological aspect this article is centred around is the following question:

\begin{ques*}
Under which conditions on $G, H$ and the inclusion is a
homogeneous space (respectively a biquotient) (non-)formal?
\end{ques*}

Recall that a topological space $X$ is called \emph{formal}, if the information contained in its rational cohomology algebra is ``the same'' as its entire \emph{rational homotopy type}. In particular, the rational homotopy groups $\pi_*(X)\otimes \qq$ of a simply-connected formal $X$ may be computed from its cohomology algebra. Let us give the precise definition:
\begin{defi*}
The commutative differential graded algebra $(A,\dif)$ is called \emph{formal} if it is weakly equivalent
to the cohomology algebra $(H(A,\qq),0)$, i.e.~if there is a chain of quasi-isomorphisms
\begin{align*}
(A,\dif) \xto{\simeq} \dots \xleftarrow{\simeq} \dots \xto{\simeq} \dots
\xleftarrow{\simeq} (H(A,\qq),0)
\end{align*}
We call a path-connected topological space \emph{formal} if its
rational homotopy type is a formal consequence of its rational
cohomology algebra, i.e.~if $(\APL(X),\dif)$, the commutative differential graded algebra of polynomial differential forms on $X$, is formal
\end{defi*}
Formality is one of the most important and most discussed topics in Rational Homotopy Theory. Conjecturally, it forms an obstruction to the existence of metrics of positive curvature on manifolds. Besides, it can be used to distinguish K\"ahler manifolds from symplectic manifolds, and it is an obstruction to \emph{geometric formality}---which also was extensively studied on homogeneous spaces.

An elaborate discussion of the formality of homogeneous spaces was
given in the classical book \cite{GHV76}. Amongst others this
resulted in long lists of formal homogeneous spaces
(cf.~\cite{GHV76}.XI, p.~492-497). It is well known that symmetric spaces of compact type are formal (cf.~\cite{FHT01}.12.3, p.~162), $N$-symmetric spaces are
formal (cf.~\cite{Ste02}, Main Theorem, p.~40, for the precise
statement, \cite{KT03}). (The first mentioned spaces are even geometrically formal; this is not true in general in the latter case.) Moreover, if $\rk G=\rk H$ the space $\biq{G}{H}$ is \emph{positively elliptic} (or \emph{$F_0$}), i.e.~rationally elliptic with positive Euler characteristic and consequently, formal.

However, it seems that only very few examples of
non-formal homogeneous spaces are known. For example consider \cite{GHV76}.XI.5, p.~486-491. (One
example was generalised to a parametrised family in \cite{Tra93}. A few further examples are cited in example \cite{KT11}.1, p.~158.)
Exemplarily, the following homogeneous spaces are known to be
non-formal for $p,q\geq 3$ and $n\geq 5$:
\begin{align*}
\frac{\SU(pq)}{\SU(p)\times\SU(q)}, \qquad \frac{\Sp(n)}{\SU(n)}
\end{align*}
In this article we shall discuss mainly three principles of how to construct non-formal homogeneous spaces respectively biquotients. In particular, this will produce uniform proofs for the non-formality of the known examples and it will result in more examples like
\begin{align*}
\frac{\SU(p+q)}{\SU(p)\times \SU(q)}, \qquad \frac{\SO(2n)}{\SU(n)}
\end{align*}
for $p+q\geq 4$, $n\geq 8$ just to mention some simple cases. We refer the reader to theorems \ref{HOMtheo01} and \ref{theo01} as well as example \ref{ex02} for several more non-formal homogeneous spaces.

\vspace{5mm}

We provide the following tools to construct these spaces.
\begin{mainprop}\label{propA}
Let $H\In G$ and $K\In H\times H$ be compact connected Lie
groups. Suppose that the inclusion of $H$ into $G$
induces an injective morphism on rational homotopy groups,
i.e.~$\pi_*(H)\otimes \qq\hto{}\pi_*(G)\otimes \qq$.

Then $\biq{G}{K}$ is formal if and only if $\biq{H}{K}$ is formal.
\end{mainprop}

For the next theorem we need to recall the following: A space $E$ possesses the \emph{Hard-Lefschetz property} if there exists a closed $2$-form $l\in H^2(E,\rr)$ such that for all $k\in \nn$
\begin{align*}
L^k: H^{n-k}(E,\rr) \to H^{n+k}(E,\rr) \qquad L^k([\alpha])=[l^k\wedge \alpha]
\end{align*}
is an isomorphism. The most prominent examples of \emph{Hard-Lefschetz manifolds}, i.e.~manifolds with the Hard-Lefschetz property, are K\"ahler manifolds.
Note that a space is formal (over $\qq$) if and only if it is so over each field extension of $\qq$ (cf.~\cite{FHT01}, p.~156 and theorem 12.1, p.~316); thus we need not worry about coefficients here.

The next theorem will provide the main source for finding non-formal homogeneous spaces.
\begin{main}\label{theoA}
Let $E^{2n+1}$ be a simply-connected space with finite dimensional rational cohomology and let
$B^{2n}$ be a simply-connected Hard-Lefschetz space (for $n\geq 1$). Suppose
there is a fibration $\s^{1}\hto{} E\xto{p} B$. Assume further the
following to hold true:
\begin{itemize}
\item The Euler class of $p$ is a non-vanishing multiple of the
K\"ahler class $l$ of $B$ in rational cohomology.
\item
The rational cohomology of $B$ is concentrated in even degrees only.
\item
The rational homotopy groups of $B$ are concentrated in degrees
\linebreak[4]smaller or equal to $n$.
\end{itemize}
Then $E$ is a non-formal space.
\end{main}
We remark that although this theorem requires extremely special conditions, it fits perfectly to the case of homogeneous spaces (and bioquotients). As we shall see, all these requirements can easily be fulfilled in a large class of cases.

This theorem has a ``mirror version'', which does give a better characterisation of formality in this situation. Recall that a space is called an \emph{$F_0$-space} if it is rationally elliptic with positive Euler characteristic.
\begin{main}\label{theoB}
Let $B$ be a positively elliptic simply-connected Hard-Lefschetz space of dimension $2n$. Let
\begin{align*}
\s^1\hto{} E \to B
\end{align*}
be a fibration with simply-connected total space $E$ of formal dimension $2n+1$. Suppose that the Euler class of the fibration equals the Hard-Lefschetz class (up to non-trivial multiples in the second rational cohomology group).

Then $E$ is formal if and only if
$E$ splits rationally as a product with one factor an odd-dimensional rational sphere of dimension greater than or equal to $n$ and with the other factor being an $F_0$-space.
\end{main}
We remark that an elliptic space is positively elliptic if and only if its rational cohomology is concentrated in even degrees---cf.~proposition \cite{FHT01}.32.10, p.~444.

The last criterion will be
\begin{mainprop}\label{propD}
Let $G$ be a rationally elliptic H-space of finite type. Suppose
it admits a free group action by a compact connected Lie group $H$.
Then $G/H$ is formal if and only if $G/T_H$ is formal.
\end{mainprop}
This implies, in particular, that $\biq{G}{H}$ is formal if and only if so is $\biq{G}{T_H}$---a generalisation of the homogeneous case (cf.~\cite{Oni94}, p.~212). This principle follows from our much more general result in \cite{AK11}; yet, we provide a much simpler proof here. (In order to avoid confusion, we remark that the term ``H-space" refers to the existence of a multiplication up to homotopy and not to the action of the Lie group $H$.)

\vspace{5mm}

It is a classical result (see \cite{NM78}) that simply-connected compact
manifolds of dimension at most six are formal spaces. In
\cite{FM04} Fern\'andez and Mu\~noz pose the question whether
there are non-formal simply-connected compact manifolds in
dimensions seven and higher. They answer this in the affirmative by
constructing seven-dimensional and eight-dimensional examples. The
requested example in a certain dimension above dimension seven is
then given by a direct product with the corresponding
even-dimensional sphere.

A next goal is to find features which usually tend to be favourable to establishing formality and to construct non-formal examples bearing these properties. One direction in this vein is to find
highly-connected examples (cf.~for example \cite{FM06}). Note that
the methods of constructing these manifolds involve surgery theory,
in particular.

We shall use a different approach towards this existence problem of ``geometrically nice'' non-formal spaces in every large dimension.

Combining the Bott conjecture and the Hopf conjecture, every manifold of positive curvature should be positively elliptic and formal, in particular. The $G$-invariant metric on homogeneous spaces $G/H$ of compact Lie groups has non-negative curvature. However, we shall see that from dimension $72$ on in each dimension there exists a compact manifold of non-negative curvature, which is not formal, but elliptic.

\begin{main}\label{theoE}
In every dimension $d\geq 72$ there is an irreducible simply-connected compact homogeneous
space which is not formal. In
particular, every such space constitutes a rationally elliptic
non-formal manifold admitting a metric of non-negative curvature.
\end{main}

\vspace{3mm}

As a reference for Rational Homotopy Theory we draw on the textbook \cite{FHT01} and we shall follow its terminology and notation.

\textbf{All the Lie groups under consideration will be compact connected. All commutative differential graded algebras---and cohomology, in particular---will be taken with rational coefficients unless stated differently.}

\vspace{3mm}

\str In section \ref{sec01} we shall prove the first two contruction principles---proposition \ref{propA}, theorems \ref{theoA} and \ref{theoB}---and we provide a large list of non-formal homogeneous examples. Section \ref{sec02} is devoted to providing the necessary examples for the proof of theorem \ref{theoE}. Finally, in section \ref{sec03} we show proposition \ref{propD} via a discussion of formality in special fibrations.

\vspace{3mm}

\ack The author is very grateful to Anand Dessai for various
fruitful discussions and to Jim Stasheff for commenting on a previous version of this article.

\vspace{3mm}


\section{Construction principles}\label{sec01}

In theory, given the groups $G$, $H$ and the inclusion
of $H$ into $G$ it is possible to compute whether the space $G/H$ is
formal or not. Consider the appendix \ref{app} for an exemplary computation.
However, understanding the topological nature of the
inclusion may be a non-trivial task. Moreover, computational complexity theory enters the stage:
In \cite{GL03} the following result was established: Given a
rationally elliptic simply-connected space, computing its Betti numbers from
its minimal Sullivan model---which
serves as an encoding of the space---is an NP-hard problem. In \cite{Ama11b} it is proved that computing the rational cup-length and the rational Lusternik--Schnirelmann category of an elliptic space is also NP-hard. So on an elliptic space of large dimension it should be rather challenging
to decide whether the space is formal or not. Besides, it is ``desirable'' (cf.~\cite{Ste02}, p.~38) to give criteria
on invariants of $G$, $H$ and the embedding $H\hto{} G$ by which a
direct identification of the homogeneous space $G/H$ as being formal
or not may be achieved. (In the article \cite{Ste02}, this was done by
identifying $N$-symmetric spaces as formal. See also \cite{KT03} for
the same result.)

Thus finding a priori arguments which allow to
identify infinite non-formal families should be worth-while the
effort. Let us begin to do so by proving proposition \ref{propA}.
\begin{proof}[\textsc{Proof of proposition \ref{propA}}]
Compact connected Lie groups are simple \linebreak[4]spaces. Choose minimal Sullivan models $(\Lambda V_G,0)$,
$(\Lambda V_H,0)$ and $(\Lambda V_K,0)$ for $G$,
$H$ and $K$. Let $x_1,\dots,x_n$ be a homogeneous basis of $V_H$.
Since the inclusion $H\hto{} G$ induces an injective morphism on
rational homotopy groups, we may choose $x_1',\dots, x_k'$ such that
$x_1,\dots, x_n,x'_1,\dots, x'_k$ is a homogeneous basis of $V_G$.
(For this we identify generators of homotopy groups with the
$x_i,x'_i$ up to duality.)

So we obtain
\begin{align*}
H^*(G\times G)\cong \Lambda \langle x_1,\dots,x_n,x'_1,\dots,x'_k
,y_1,\dots,y_n,y'_1,\dots y'_k\rangle
\end{align*}
(The $y_i, y'_i$ are constructed just like the $x_i, x'_i$, i.e.~as
a ``formal copy''.)

Hence a model for the biquotient $\biq{G}{K}$ is given by
(cf.~proposition \cite{Kap}.1, p.~2)
\begin{align*}
(\Lambda V_K^{+1}\otimes \Lambda \langle
q_1,\dots,q_n,q'_1,\dots,q'_k \rangle,\dif)
\end{align*}
(The degrees of $V_K$ are shifted by $+1$.)

Denote the inclusion $K\hto{i_1}H\times H \hto{i_2} G\times G$ by
$\phi=i_2\circ i_1$. The differential $\dif$ on the $q_i$ is given
by $\dif(q_i)=H^*(\mathbf{B}\phi)(\bar x_i-\bar y_i)$, where $\bar
x_i,\bar y_i\in H^*(\mathbf{B}G)$ is the class corresponding to
$x_i,y_i\in H^*(G)$. (The analogue holds for the $q'_i$.) By
functoriality the morphism $H^*(\mathbf{B}\phi)$ factors through
\begin{align*}
H^*(\mathbf{B}G\times \mathbf{B}G)\xto{H^*(\textbf{B}i_2)}
H^*(\mathbf{B}H\times \mathbf{B}H)\xto{H^*(\textbf{B}i_1)}
H^*(\mathbf{B}K)
\end{align*}
In particular, $H^*(\mathbf{B}\phi)=H^*(\mathbf{B}i_1)\circ H^*(\mathbf{B}i_2)$ and
\begin{align*}
\dif(q'_i)=H^*(\mathbf{B}i_1)(H^*(\mathbf{B}i_2)(\bar
x'_i-\bar y'_i))
\end{align*}
for $1\leq i\leq k$. The inclusion $i_2$ is the
product $i_2=i_2|_H\times i_2|_H$  by definition. Thus we obtain a decomposition
$H^*(\textbf{B}i_2)=H^*(\textbf{B}(i_2|_H))\otimes
H^*(\textbf{B}(i_2|_H))$ and $H^*(\mathbf{B}i_2)(\bar x'_i)\in
\Lambda \langle x_1,\dots,x_n\rangle$; i.e.~it is a linear
combination in products of the $x_i$. By the splitting of
$H^*(\textbf{B}i_2)$ we obtain that $H^*(\mathbf{B}i_2)(\bar
y'_i)\in \Lambda \langle y_1,\dots,y_n\rangle$ is the same
linear combination with the $x_i$ replaced by the $y_i$. This
implies that
\begin{align*}
\dif(q'_i)=H^*(\mathbf{B}i_1)(H^*(\mathbf{B}i_2)(\bar
x'_i-\bar y'_i))\in \im H^*(\mathbf{B}i_1)=
\dif(\Lambda \langle q_1, \dots,q_n\rangle)
\end{align*}
(In fact, the very same linear combination as above with the $x_i$
now replaced by the $q_i$ will serve as a preimage of $\dif(q_i')$
under $\dif|_{\Lambda \langle q_1,\dots,q_n\rangle}$.)

Thus for each $1\leq i\leq k$ there is an element $z_i\in \Lambda
\langle q_1, \dots,q_n\rangle$ with the property that
$\dif(q'_i-z_i)=0$. Set $\tilde q_i:=q'_i-z_i$. We have an
isomorphism of commutative differential graded algebras
\begin{align*}
\sigma:& (\Lambda V_K^{+1}\otimes \Lambda \langle
q_1,\dots,q_n,q'_1,\dots,q'_k \rangle,\dif)\\&\xto{\cong} (\Lambda
V_K^{+1}\otimes \Lambda \langle q_1,\dots,q_n,\tilde
q_1,\dots,\tilde q_k \rangle,\dif)
\end{align*}
induced by the ``identity'' $q_i'\mapsto \tilde q_i+ z_i$ for all
$1\leq i\leq k$. (By abuse of notation we now consider the $\tilde
q_i$ with $\dif \tilde q_i:=\dif (q'_i-z_i)=0$ as abstract elements
in the graded vector space upon which the algebra is built.) This
morphism is an isomorphism of commutative graded algebras which
commutes with differentials:
\begin{align*}
\dif(\sigma(q'_i))= \dif(\tilde q_i+
z_i)=\dif(q'_i)=\sigma(\dif(q'_i))
\end{align*}
as $\sigma|_{\Lambda \langle q_1,\dots, q_k\rangle}=\id$.

Thus we obtain a
quasi-isomorphism
\begin{align}\label{HOMeqn09}
\nonumber\APL(\biq{G}{K})  \simeq &(\Lambda
V_K^{+1}\otimes \Lambda \langle q_1,\dots,q_n,q'_1,\dots,q'_k
\rangle,\dif)\\\cong & (\Lambda V_K^{+1}\otimes \Lambda
\langle
q_1,\dots,q_n,\tilde q_1,\dots,\tilde q_k \rangle,\dif)\\
\nonumber= & (\Lambda V_K^{+1} \otimes \Lambda \langle
q_1,\dots,q_n\rangle,\dif) \otimes (\Lambda\langle\tilde
q_1,\dots,\tilde q_k \rangle,0)
\end{align}
The algebra $(\Lambda V_K^{+1} \otimes \Lambda \langle
q_1,\dots,q_n\rangle,\dif)$ is a model for $\biq{H}{K}$, since
\linebreak[4] $x_1,\dots, x_n, y_1,\dots,y_n$ is a basis of $V_H$ and since its
differential $\dif$ corresponds to $H^*(\textbf{B}i_1)$. Hence the
last algebra in \eqref{HOMeqn09} is rationally the product of a
model of $\biq{H}{K}$ and a formal algebra. Thus it is formal if and
only if $\biq{H}{K}$ is formal (cf.~theorem \cite{Ama09}.1.47, p.~39).
\end{proof}
\begin{cor}\label{HOMcor01}
Let $K$ be a compact Lie subgroup of the Lie group $G$, which itself
is a Lie subgroup of a Lie group $\tilde G$. Table \ref{table01}
gives pairs of groups $G$ and $\tilde G$ together with relevant
relations and the type of the inclusion such that it holds: The
homogenenous space $G/K$ is formal if and only if $\tilde
G/K$ is formal.
\begin{table}
\centering \caption{Inclusions of Lie groups}
\label{table01}
\begin{tabular}{ c@{\hspace{1mm}} |@{\hspace{1mm}}c| @{\hspace{1mm}}l|@{\hspace{1mm}}c}
&&& \\[-3.5mm]$G$ & $\tilde G$ & with & type of embedding\\
\hline
 &&&\\[-3mm]  $\SO(n)$ & $\SO(N)$& $n\geq 1$, $n$ odd, $N\geq n$, $N$ odd & blockwise\\
&&&\\[-3mm]$\SO(n)$ & $\SO(N)$& $n\geq 1$, $n$ odd, $N\geq n+1$, $N$ even & blockwise\\
&&&\\[-3mm]$\SU(n)$ & $\SU(N)$ & $n>1$, $N\geq n$& blockwise\\
&&&\\[-3mm]$\Sp(n)$ & $\Sp(N)$ & $n\geq 1$, $N\geq n$& blockwise \\
 &&&\\[-3mm]$\SO(n)$ & $\SU(n)$ &  $n\geq 1$, $n$ odd &induced
\\&&&componentwise
\\&&&by $\rr\hto{}\cc$
\\
 &&&\\[-3mm]$\Sp(n)$ & $\SU(2n)$ &  $n\geq 1$ & induced by the
 \\&&& identification \\&&&$\hh\hto{} \cc^{2\times 2}$\\
\end{tabular}
\end{table}
\end{cor}
\begin{prf}
For the minimal models of the relevant Lie groups
and further reasoning on their inclusions see \cite{FHT01}.15, p.~220.
Lie groups are formal and so are the depicted
inclusions.
More precisely,
we may identify minimal models with
their cohomology algebra and induced maps on minimal models with
induced maps in cohomology. Thus it suffices to see that the given
inclusions induce surjective morphisms in cohomology.

The blockwise inclusions are surjective in cohomology due to theorem
\cite{MT91}.6.5.(4), p.~148.

The inclusion $\SO(n)\hto{} \SU(n)$ induces a surjective morphism by
\linebreak[4]\cite{MT91}.6.7.(2), p.~149. So does the inclusion $\Sp(n)\hto{}
\SU(2n)$ by \cite{MT91}.6.7.(1), p.~149.
\end{prf}
We remark that the chain of inclusions $\SO(n)\In\U(n)\In\Sp(n)$ does \emph{not} induce a surjective morphism in cohomology
(cf.~\cite{MT91}.5.8.(1), p.~138 and \cite{MT91}.6.7.(2), p.~149).
Neither does the chain $\Sp(n)\In \U(2n)\In \SO(4n)$ nor the chain
$\Sp(n)\In \U(2n)\In \SO(4n+1)$ in general by a reasoning taking
into account theorem \cite{MT91}.6.11, p.~153, which lets us conclude
that the transgression in the rationalised long
exact homotopy sequence for $\U(n)\hto{} \SO(2n) \to \SO(2n)/\U(n)$
is surjective. Let us give a first example, which will be largely generalised later.
\begin{ex}\label{ex01}
\begin{itemize}
\item
Using the computation from the appendix \ref{app} and proposition \ref{propA} we
see that the space
\begin{align*}
\frac{\SU(n)}{\SU(3)\times \SU(3)}
\end{align*}
is non-formal for $n\geq 6$.
\item
We also obtain pretty simple proofs for formality: For example, the
spaces
\begin{align*}
\frac{\U(n)}{\U(k_1)\times \dots \times \U(k_l)} \qquad
\textrm{and}\qquad \frac{\SO(2n+1)}{\SO(2k_1)\times \dots \times
\SO(2k_l)}
\end{align*}
with $\sum_{i=1}^l k_i\leq n$ are formal: In the equal rank case
$\sum_{i=1}^l k_i=n'$ the spaces are positively elliptic and formality holds. Then apply proposition \ref{propA}
to the inclusions $\U(n')\hto{}\U(n)$ and
$\SO(2n'+1)\hto{}\SO(2n+1)$.
\end{itemize}
\end{ex}

Let us now prove theorem \ref{theoA} which serves to derive
non-formality via additional topological structure.
This will simplify proofs for some of the known
non-formal examples whilst producing a whole variety of further
ones. We cite the following technical lemma from \cite{DGMS75}, theorem 4.1, p.~261, and from lemma \cite{FM05}.2.7 p.~154.
\begin{lemma}\label{lemma01}
A minimal model $(\Lambda V,\dif)$
is formal if and only if it can be written (up to isomorphism) such that there is in each $V^i$ a complement $N^i$
to the subspace of $\dif$-closed elements $C^i$ with $V^i= C^i\oplus
N^i$ and such that any closed form in the ideal $I(\bigoplus N^i)$
generated by $\bigoplus N^i$ in $\Lambda V$ is exact.
\end{lemma}
\vproof

\begin{proof}[\textsc{Proof of Theorem \ref{theoA}}]
Since the K\"ahler class and the Euler class rationally are
non-trivial multiples, we may formulate the Gysin sequence with
rational coefficients for $p$ as follows:
\begin{align*}
\dots \to H^p(B) \xto{\cup [l]} H^{p+2}(B) \to H^{p+2}(E)\to
H^{p+1}(B)\xto{\cup [l]}\dots
\end{align*}
(The fibration is oriented, since $B$ is simply-connected.) The
Hard-Lefschetz property of $B$ implies that taking the cup-product with $[l]$ is
injective in degrees $p<n$. So the sequence splits, yielding
\begin{align*}
H^{p+2}(B)&=H^p(B)\oplus H^{p+2}(E)
\end{align*}
for $-2\leq p\leq n-2$. Since we assumed the odd-dimensional rational cohomology groups of $B$ to vanish, we obtain in particular that
\begin{align}
\label{HOMeqn05}
H^p(E)&=0 &\textrm{for odd } p\leq n
\intertext{Let $(\Lambda V_E,\dif_E)$ be a minimal model of $E$. As $E$ was supposed to be simply-connected, we clearly have $\pi_p(E)\otimes \qq\cong V_E^p$ up to duality. By the long exact sequence of homotopy groups we obtain}
\nonumber \pi_p(E)&\cong\pi_p(B) &\textrm{for } p\geq 3
\intertext{As we assumed}
\nonumber\pi_p(B)\otimes \qq&=0 &\textrm{for }p>n
\intertext{we see that $V_E^{>n}=0$, if
$n\geq 2$. If $n=1$, we clearly still have $V_E^{\geq 3}$=0. So for arbitrary $n\geq 1$ we obtain}
\label{HOMeqn02}
V_E^p&=0 &\textrm{for odd }p>n
\end{align}

\vspace{5mm}

Assume now that $E$ is formal. This will lead to a
contradiction. By lemma \ref{lemma01} we may split (a suitable minimal model $(\Lambda V_E,\dif_E)$ of $E$ as)
$V_E=C_E\oplus N_E$ with $C_E=\ker \dif_E|_{V_E}$ and with the
property that every closed element in $I(N_E)$ is exact.

By observation \eqref{HOMeqn05} we know that every closed element of
odd degree in $(\Lambda V_E)^{\leq n}$ is exact. Thus by
the minimality of the model we directly derive that
\begin{align}
C_E^{\leq n}&=(C_E^{\leq n})^\textrm{even}
\intertext{Together with \eqref{HOMeqn02} this implies that}
\label{HOMeqn03} C_E&=C_E^\textrm{even} \intertext{is concentrated
in even degrees only. Hence so is $\Lambda C_E$. Thus, using the
splitting}
\nonumber\Lambda V_E&=\Lambda C_E \oplus I(N_E)
\intertext{we derive that}
\label{HOMeqn04}
H^\textrm{odd}(E)&=\bigg(\frac{\ker \dif_E|_{I(N)}}{\im \dif_E}\bigg)^\textrm{odd}=0\
\end{align}
from \eqref{HOMeqn03} and from the formality of $E$.

However, by assumption, the space $E$ is an odd-dimensional simply-connected rationally elliptic space
Thus its rational cohomology satisfies Poincar\'e duality by theorem \cite{FHT91}.A, p.~70. The space $E$ possesses a volume form which hence generates
$H^{2n+1}(E)\cong\qq\neq 0$. This contradicts formula \eqref{HOMeqn04}.
Thus $E$ is non-formal.
\end{proof}
\begin{rem}\label{HOMrem03}
\begin{itemize}
\item
Actually, we do not only prove the non-formality of the total space; we can even derive the existence of a non-vanishing Massey product of odd degree. (Note that there are examples of non-formal spaces with vanishing Massey products.)
\item
We see that the proof works equally well in the following setting:
Assume the spaces $E$ and $B$ not to be simply-connected but to be compact connected manifolds with a finite fundamental group and to be \emph{simple spaces}, i.e.~the action of the fundamental group on all
homotopy groups (in positive degree) is trivial. Indeed, the spaces
then are orientable and satisfy Poincar\'e duality. Since we may choose the Euler class of the $\s^1$-fibration up to (non-trivial) multiples, we then may assume that the fibration is orientable (up to coverings of the fibre), and the Gysin sequence can be applied.
Moreover, Rational Homotopy Theory is equally applicable to simple spaces.

From proposition \cite{Oni94}.1.17, p.~84, we cite that the homogeneous
space $G/H$ of a compact connected Lie group $G$ with a connected
closed Lie subgroup $H$ is simple.
\item
K\"ahler manifolds are clearly not the only source of examples of $B$. For example, the lemma applies to the case of certain Donaldson submanifolds
(cf.~\cite{FM05}) or biquotients (cf.~\cite{Kap}).

We may relax the assumptions of the lemma even further: We need
not require the rational homotopy groups of $B$ to be concentrated
in degrees smaller or equal to $n$ as long as the ones concentrated in
odd degrees above $n$ correspond to relations in the minimal model
of $E$ via the long exact homotopy sequence; i.e.~they are not in
the image of the dual of the rationalised Hurewicz
homomorphism.
Then the proof does not undergo any severe
adaptation: We just obtain $C_E^p=0$ instead of $V_E^p=0$ for odd
$p>n$. This extends the lemma to not necessarily elliptic spaces $E$
and $B$.
\item
In theory, the approach of this theorem is not restricted to $\s^1$-fibrations. So, for
example, using the Kraines form in degree four instead of the
K\"ahler form in degree two, one may formulate a version for Positive
Quaternion K\"ahler Manifolds---however, this seems to be less fruitful.
\end{itemize}
\end{rem}
See the corollary on \cite{Oni94}.13, p.~221 and the examples below
for a related result in the category of homogeneous spaces.

Via theorem \ref{theoB} we cast a little more light on the situation depicted in theorem \ref{theoA} and we characterise when formality occurs. So before we begin to use theorem \ref{theoA} in order to find several examples of non-formal homogeneous spaces---which clearly provides the motivation for  the theorem---let us prove theorem \ref{theoB}. We recall the following
\begin{lemma}\label{lemma02}
A pure elliptic minimal Sullivan
algebra $(\Lambda V,\dif)$ is
formal if and
only if it is of the form
\begin{align}\label{PUReqn03}
(\Lambda V,\dif)\cong  (\Lambda V',\dif) \otimes
(\Lambda \langle z_1,\dots, z_l\rangle,0)
\end{align}
(for maximal such $l$) with a pure minimal Sullivan algebra
$(\Lambda V',\dif)$ of positive Euler characteristic---which is
automatically formal then---and with odd degree generators $z_i$.
\end{lemma}
\begin{prf}
Choose a basis $z_1,\dots, z_m$ of $V^\textrm{odd}$ with the
property that $\dif z_i=0$ for \linebreak[4]$1\leq i \leq l$---and some fixed
$1\leq l\leq m$---and that $\dif|_{\langle
z_{l+1},\dots,z_m\rangle}$ is injective. From pureness it follows:
\begin{align*}
\dif(\langle z_{1},\dots,z_l\rangle)&=0\\
\dif(\langle z_{l+1},\dots,z_m\rangle)&\in \Lambda V^\textrm{even}\\
\dif(V^\textrm{even})&=0
\end{align*}
Thus we obtain
\begin{align*}
(\Lambda V,\dif)=  (\Lambda (V^\textrm{even}\oplus\langle z_{l+1},\dots,z_m
\rangle),\dif) \otimes (\Lambda \langle z_{1},\dots,z_l
\rangle,0)
\end{align*}
Set $V':=V^\textrm{even}\oplus\langle z_{l+1},\dots,z_m
\rangle$. The minimality of $(\Lambda V,\dif)$ enforces
the minimality of $(\Lambda V',\dif)$ and $(\Lambda
V',\dif)$ is again a pure minimal Sullivan algebra.

By theorem \cite{Ama09}.1.47, p.~39, the formality of $(\Lambda V,\dif)$ now is equivalent
to the formality of $(\Lambda V',\dif)$.

If $(\Lambda V',\dif)$ has positive Euler characteristic, it is positively elliptic and formal---cf.~theorem \cite{Ama09}.1.52, p.~43.
Thus $(\Lambda V,\dif)$ then is formal, too.

For the reverse implication we use contraposition: Suppose the Euler
characteristic of $(\Lambda V',\dif)$ to vanish. We shall
show that the algebra is non-formal. The algebra $(V',\dif)$
has the property that $\dif ((V')^\textrm{even})=0$ and that
$\dif|_{(V')^\textrm{odd}}$ is injective. Thus, setting
$C:=(V')^\textrm{even}$ and $N:=(V')^\textrm{odd}$, yields a
decomposition $V'=C\oplus N$ as required in lemma
\ref{lemma01}; this decomposition is uniquely determined up to isomorphism.

Due to
$\chi(\Lambda V',\dif)=0$ there is a closed and non-exact
element $x$ in $(\Lambda V',\dif)$ that has odd degree.
Therefore $x$ necessarily lies in $I(N)$. Thus $(\Lambda V',\dif)$ is not formal and neither is $(\Lambda V,\dif)$.
\end{prf}

\begin{proof}[\textsc{Proof of theorem \ref{theoB}}]
Choose minimal models $(\Lambda V,\dif)$ for $B$ and \linebreak[4]$(\Lambda \langle s\rangle, 0)$ for $\s^1$. Then we obtain a model for the fibration---cf.~proposition \cite{FHT01}.15.5, p.~198:
\begin{align*}
(\Lambda V  \otimes \Lambda \langle s\rangle ,\dif)\xto{\simeq} (\APL(E),\dif)
\end{align*}
with $\dif s=l\in V^2$ and $[l]\in H^2(Z)$ representing the Hard-Lefschetz class.

Obviously, the model of the fibration is not minimal. However, we may construct its minimal model. For this we write the model as
\begin{align*}
(\Lambda V' \otimes \Lambda \langle s,l\rangle ,\dif)
\end{align*}
with $V'\In V$ a homogeneous complement of $\langle l\rangle$. This point of view reveals this model as the model of a rational fibration with fibre $\Lambda V'$ and contractible base space $(\Lambda \langle s, l\rangle,\dif')$. Thus we obtain (quasi-)isomorphisms
\begin{align*}
&(\Lambda V \otimes \Lambda \langle s\rangle ,\dif)\\
\cong & (\Lambda V' \otimes \Lambda \langle s,l\rangle ,\dif)
\\\simeq & (\Lambda V' \otimes \Lambda \langle s,l\rangle ,\dif)/(\Lambda \langle s, l\rangle,\dif')
\\\simeq & (\Lambda V',\bar \dif)
\end{align*}
where $\bar \dif$ is the projection of the differential $\dif$ to $V'$. Since $(\Lambda V,\dif)$ is minimal, so is $(\Lambda V',\bar \dif)$. Hence we have found a minimal model of $E$.

Applying the Gysin sequence to the fibration
\begin{align*}
\ldots\to H^k(E)\to H^{k-1}(B)\xto{\cup [l]}H^{k+1}(B)\to H^{k+1}(E)\to \ldots
\end{align*}
yields
\begin{align}\label{eqn04}
H^k(E)=0
\end{align}
for odd $k<n$. This follows from the assumption that $B$ is positively elliptic, i.e.~that its rational cohomology is concentrated in even degrees. Moreover, we use the assumption that the Euler class of the fibration is a rational multiple of the Hard-Lefschetz class $[l]$. Due to Hard-Lefschetz, cupping with $[l]$ is injective in degrees smaller than $n$.

We may now characterise formality. The model $(\Lambda V',\bar\dif)$ has odd formal dimension $2n+1$. Since the space $E$ is elliptic by the long exact homotopy sequence of the fibration, it satisfies Poincar\'e duality. In particular, there is a non-vanishing cohomology class of odd degree.

Since the space $B$ is positively elliptic, it is pure. Consequently, also the model $(\Lambda V',\bar \dif)$ is pure. By lemma \ref{lemma02} and the long exact sequence of the fibration we obtain the following: The formality of $E$ is equivalent to the existence of a decomposition
\begin{align*}
(\Lambda V',\bar\dif)\cong (\Lambda V'',\tilde\dif) \otimes (\Lambda \langle \tilde n\rangle, 0)
\end{align*}
for  $0\neq \tilde n\in (V')^\textrm{odd}$ and a homogeneous complement $V''$ of $\langle \tilde n\rangle$ in $V'$. Thus $E$ is formal if and only if its minimal model splits as a product of a formal space $(\Lambda V'',\tilde\dif)$ and the model of an odd sphere.

In this case we observe that
\begin{align*}
\dim  (V'')^\textrm{odd}=\dim (V'')^\textrm{even}
\end{align*}
Indeed, this was true for $V$ already, since the base is positively elliptic. The equality above then follows, since $V''$ results from $V$ by splitting of a class of even degree, $l$, and a class of odd degree, $\tilde n$. Thus, the base $B$ is again positively elliptic (cf.~proposition \cite{FHT01}.32.10, p.~444), since the
cohomology $H(\Lambda V'',\tilde\dif)$ is finite-dimensional by construction. Hence $(\Lambda V'',\tilde\dif)$ is an $F_0$-space and consequently formal.

We conclude that this shows that $E$ is formal if and only if it splits as a product of formal spaces, which then are necessarily an $F_0$-space and a rational sphere of odd dimension at least $n$---cf.~equation \eqref{eqn04}. This proves the result.
\end{proof}

A main step in the construction of our non-formal homogeneous spaces will be to consider certain $\s^1$-bundles over homogeneous Hermitian spaces. For this we cite the following simply-connected compact homogeneous
K\"ahler manifolds from \cite{Bes08}.8.H, p.~229--234:
\begin{align}\label{HOMeqn07}
\nonumber&\frac{\Sp(n)}{\U(p_1)\times \dots \times \U(p_q)\times \Sp(l)} & \textrm{for } \sum_{i=1}^q p_i+l=n\\
\nonumber&\frac{\SO(2n+1)}{\U(p_1)\times \dots \times \U(p_q)\times \SO(2l+1)}& \textrm{for } \sum_{i=1}^q p_i+l=n\\
&\frac{\SO(2n)}{\U(p_1)\times \dots \times \U(p_q)\times \SO(2l)}& \textrm{for } \sum_{i=1}^q p_i+l=n\\
\nonumber&\frac{\SO(2n)}{\U(p_1)\times \dots \times \U(p_{q-1})\times \tilde\U(p_q)}& \textrm{for } \sum_{i=1}^q p_i=n\\
\nonumber&\frac{\SU(n)}{\mathbf{S}(\U(p_1)\times \dots \times
\U(p_q))}& \textrm{for } \sum_{i=1}^q p_i=n
\end{align}
where $\tilde \U(p_q)\In \SO(2p_q)$ is the unitary group with
respect to a slightly altered complex structure of $\rr^{2p_q}$
(cf.~\cite{Bes08}.8.113, p.~231).

One important observation will be that the dimension of these spaces grows quadratically in $n$, whereas the top dimension of their rational homotopy groups grows linearly. This enables us to apply theorem \ref{theoA}.

In table \ref{HOMtable03} we shall use the convention that
$\sum_{i=1}^t k_i=p$, $\sum_i l_i=l$ and $p_1+p_2=p$. Moreover, we also require
$p,p_1,p_2>0$ in the table.
\begin{theo}\label{HOMtheo01}
The homogeneous spaces in table \ref{HOMtable03} are non-formal.
\begin{table}[h]
\centering \caption{Examples of non-formal homogeneous spaces}
\label{HOMtable03}
\begin{tabular}{ c@{\hspace{1mm}} |@{\hspace{1mm}}l}
 homogeneous space & with \\
\hline  &\\$\frac{\Sp(N)}{\mathbf{S}(\U(k_1)\times \dots
\times \U(k_t))\times \Sp(l_1)\times\dots\times \Sp(l_r) \times
\U(l_{r+1})\times \dots \times \U(l_s)}$&\tiny $p+l\leq N$, $r\geq
0$, $s\geq 0$
\\[-1mm]&\tiny$\frac{1}{2} p^2 + (2l-\frac{7}{2})p-4l+1\geq 0$ \\
&\\[-1mm]
$\frac{\SU(N)}{\mathbf{S}(\U(k_1)\times \dots \times
\U(k_t))\times \Sp(l_1)\times\dots\times \Sp(l_r) \times
\U(l_{r+1})\times \dots \times \U(l_s)}$&\tiny $2(p+l)\leq N$,
$r\geq 0$, $s\geq 0$
\\[-1mm]&\tiny$\frac{1}{2} p^2 + (2l-\frac{7}{2})p-4l+1\geq 0$ \\
&\\[-1mm]
$\frac{\SO(N)}{\mathbf{S}(\U(k_1)\times \dots \times \U(k_t)) \times
\SO(2l_1)\times \dots \times \SO(2l_r)\times \SO(2l_{r+1}+1)}$&
\tiny $2(p+l)+1\leq N$, $r+1\geq 0$\\[-1mm]
&\tiny $\frac{1}{2} p^2 + (2l-\frac{7}{2})p-4l+1\geq 0$
\\
&\\[-1mm]
$\frac{\SU(N)}{\mathbf{S}(\U(k_1)\times \dots \times \U(k_t))\times
\SO(2l_1)\times \dots \times \SO(2l_r)\times \SO(2l_{r+1}+1)}$&
\tiny $p+l=n$, $N\geq 2n+1$,\\[-1mm]&\tiny $r+1\geq 0$\\[-1mm]&\tiny
$\frac{1}{2} p^2 + (2l-\frac{7}{2})p-4l+1\geq 0$
\\
&\\[-1mm]
$\frac{\SO(2n)}{\mathbf{S}(\U(k_1)\times \dots \times \U(k_t))\times
\SO(2l_1)\times \dots \times \SO(2l_r)}$& \tiny $p+l=n$, $r\geq 0$, $n\geq 2$, $l\neq 1$ \\[-1mm]
&\tiny $\frac{1}{2}p^2+(2l-\frac{9}{2})p-4l+5\geq 0$
\\
&\\[-1mm]
$\frac{\SO(2n)}{\mathbf{S}(\tilde\U(k_1)\times \dots \times
\tilde\U(k_t))}$& \tiny $p=n$, $r\geq 0$, $n\geq 2$ \\[-1mm]
&\tiny $\frac{1}{2}n^2-\frac{9}{2}n+5\geq 0$
\\
&\\[-1mm]
$\frac{\SU(N)}{\mathbf{S}(\U(k_1)\times\dots \times \U(k_s))\times
\mathbf{S}(\U(k_{s+1})\times \dots \times \U(k_t))}$&\tiny $p_1+p_2=n\leq N$, $p_1, p_2>0$ \\[-1mm]
&\tiny $N\geq 2$, \\[-1mm]
&\tiny $p_1p_2-2(p_1+p_2)+1\geq 0$
\\
\end{tabular}
\end{table}

\end{theo}
\begin{prf}
We shall consider a certain subclass of the spaces given in
\eqref{HOMeqn07}. Their top rational homotopy group can easily be
computed using the long exact homotopy sequence of the fibration
formed by denominator, numerator and quotient. For this we shall use
that the top rational homotopy of $\Sp(n)$, $\SU(n)$, $\SO(2n)$,
$\SO(2n+1)$ lies in degree $4n-1$, $2n-1$, $4n-5$ for $n\geq 2$ and
$4n-1$ respectively. So we are able to give the spaces under consideration with their
dimensions and with the largest degree of a non-vanishing rational
homotopy group in table \ref{HOMtable04}---as always $p,p_1,p_2>0$.

\begin{table}[h]
\centering \caption{Certain Hermitian homogeneous spaces}
\label{HOMtable04}
\begin{tabular}{@{\hspace{1mm}}c@{\hspace{1mm}} @{\hspace{1mm}}r@{\hspace{1mm}}|@{\hspace{1mm}}c@{\hspace{1mm}}|@{\hspace{1mm}}c@{\hspace{1mm}}}
homogeneous space & for\quad~~\ & dimension & top rational homotopy\\
\hline  &&&\\[-1mm]$\frac{\Sp(n)}{\U(p)\times \Sp(l)}$&$p+l=n$
&$p^2+4\cdot l\cdot p+p$ &$4n-1$
\\&&&\\[-1mm] $\frac{\SO(2n+1)}{\U(p)\times \SO(2l+1)}$& $p+l=n$& $p^2+4\cdot l\cdot p+p$& $4n-1$\\
&&&\\[-1mm] $\frac{\SO(2n)}{\U(p)\times \SO(2l)}$&$ p+l=n$ & $p^2+4\cdot l\cdot p-p$& $4n-5$ for $n\geq 2$\\
&&&\\[-1mm] $\frac{\SO(2n)}{\tilde\U(n)}$ &
&$n^2-n$ &$4n-5$ for $n\geq 2$\\
&&&\\[-1mm] $\frac{\SU(n)}{\mathbf{S}(\U(p_1)\times \U(p_2))}$ &
$p_1+p_2=n$ & $2\cdot p_1\cdot p_2$&$2n-1$\\&&&\\
\end{tabular}
\end{table}

Theses homogeneous spaces  have the property that numerator and
denominator form an equal rank pair, whence they are positively elliptic and rational cohomology is concentrated in even degrees only.

All these manifolds are compact homogeneous K\"ahler manifolds by
\cite{Bes08}.8.H. Since for the following arguments the spaces one to
four behave similarly, we shall do showcase computations for the
cases one and five.

The fibre bundle
\begin{align*}
\SU(p)\hto{} \U(p) \xto{\det} \s^1
\end{align*}
lets us conclude that $\U(p)=\s^1\cdot \SU(p)=\s^1\times_{\zz_p}
\SU(p)$, where $\zz_p$ is the (multiplicative) cyclic group of $p$-th
roots of unity acting on each factor by (left) multiplication. So
the canonical projection yields a fibre bundle
\begin{align*}
\s^1/\zz_p\hto{} \frac{\Sp(n)}{(\zz_p\times_{\zz_p}\SU(p))\times
\Sp(l)} &\to \frac{\Sp(n)}{(\s_1\times_{\zz_p} \SU(p))\times \Sp(l)}
\intertext{of homogeneous spaces, which clearly is no other than}
\s^1/\zz_p\hto{} \frac{\Sp(n)}{\SU(p)\times \Sp(l)} &\to
\frac{\Sp(n)}{\U(p)\times \Sp(l)}
\end{align*}
since $\zz_p\times_{\zz_p}\SU(p)\xto{\cong}\SU(p)$---given by left
multiplication of $\zz_p$ on $\SU(p)$---is an isomorphism.

In case five we obtain a homomorphism
\begin{align*}
\s^1\times \SU(p_1)\times \SU(p_2) &\to \mathbf{S}(\U(p_1)\times
\U(p_2))\\
(x,A,B)&\mapsto (x^{p_2}\cdot A, x^{-p_1}\cdot B)
\end{align*}
with kernel formed by all the elements of the form $(a,a^{-p_2}I_{p_1},a^{p_1}I_{p_2})$ with $a^{p_1p_2}=1$. Thus the kernel is isomorphic to $\zz_{p_1p_2}$. (It is consequent
to use the terminology that $\zz_1=1$ be the trivial group.) In particular, this
defines an action of $\zz_{p_1p_2}\In\s^1$ on $\SU(p_1)\times
\SU(p_2)$ and we obtain
\begin{align*}
\mathbf{S}(\U(p_1)\times \U(p_2))=\s^1\times_{\zz_{p_1p_2}}
(\SU(p_1)\times \SU(p_2))
\end{align*}
This leads to the fibre bundle

\begin{align*}
  \s^1/\zz_{p_1p_2}\hto{}
\mbox{\fontsize{15}{19}\selectfont $\frac{\SU(n)}{\zz_{p_1p_2}\times_{\zz_{p_1p_2}}(\SU(p_1)\times
\SU(p_2))}$}  &\to \mbox{\fontsize{15}{19}\selectfont $\frac{\SU(n)}{\s_1\times_{\zz_{p_1p_2}}(
\SU(p_1)\times \SU(p_2))} $ }
\intertext{of homogeneous spaces, which clearly is no other than}
\s^1/\zz_{p_1p_2}\hto{} \frac{\SU(n)}{\SU(p_1)\times
\SU(p_2)} &\to \frac{\SU(n)}{\mathbf{S}(\U(p_1)\times \U(p_2))}
\end{align*}

\vspace{5mm}

Thus, in each case identifying $\s^1/\zz_{p}$ (for each respective $p$) with its finite covering $\s^1$, we
obtain the following sphere bundles:
{\begin{align*}
\s^1\hto{}\frac{\Sp(n)}{\SU(p)\times \Sp(l)}&\to\frac{\Sp(n)}{\U(p)\times \Sp(l)}\\
\s^1\hto{}\frac{\SO(2n+1)}{\SU(p)\times \SO(2l+1)}&\to\frac{\SO(2n+1)}{\U(p)\times  \SO(2l+1)}\\
\s^1\hto{}\frac{\SO(2n)}{\SU(p)\times \SO(2l)}&\to\frac{\SO(2n)}{\U(p)\times \SO(2l)}\\
\s^1\hto{}\frac{\SO(2n)}{\mathbf{S}\tilde\U(n)}&\to\frac{\SO(2n)}{\tilde\U(n)}\\
\s^1\hto{}\frac{\SU(n)}{\SU(p_1)\times
\SU(p_2)}&\to\frac{\SU(n)}{\mathbf{S}(\U(p_1)\times \U(p_2))}
\end{align*}}

Due to the long exact homotopy sequence, we see that the total space
$E$ of each bundle has a finite fundamental group, namely $0$ or
$\zz_2$. Thus the Euler class of each bundle does not vanish, since
otherwise the bundle would be rationally trivial. This would imply
$E\simeq_\qq B\times \s^1$---with the respective base space
$B$---and $\pi_1(E)\otimes \qq\neq 0$.

The long exact homotopy sequence associated to $H\hto{} G \to
G/H$---where $G$ is the numerator and $H$ is the denominator group
of the base space $B=G/H$ in each respective case---lets us conclude
that
\begin{align*}
H^2(B)\cong\pi_2(G/H)\otimes \qq=\qq
\end{align*}
for $l\neq 1$ and $n\geq 2$
in case three, for $n\geq 2$ in case four and without further
restrictions in the other cases. (In each respective case we then
have $\pi_2(G)\otimes \qq=0$, $\pi_1(H)\otimes \qq=\qq$ and
$\pi_1(G)\otimes \qq=0$.)

This implies that the K\"ahler class and the Euler class---both
contained in $H^2(B)$---are non-zero multiples in rational
cohomology.

In order to apply theorem \ref{theoA} we need to demand that the
top rational homotopy does not lie above half the dimension of the
space. This leads to the following restrictions in the respective
cases:
\begin{align*}
\frac{1}{2} p^2 + \bigg(2l-\frac{7}{2}\bigg)p-4l+1&\geq 0\\
\frac{1}{2} p^2 + \bigg(2l-\frac{7}{2}\bigg)p-4l+1&\geq 0\\
\frac{1}{2}p^2+\bigg(2l-\frac{9}{2}\bigg)p-4l+5&\geq 0\\
\frac{1}{2} n^2-\frac{9}{2} n +5&\geq 0\\
p_1p_2-2(p_1+p_2)+1&\geq 0
\end{align*}

Due to the first point in remark \ref{HOMrem03} we may now apply theorem
\ref{theoA} to the depicted fibre bundles. This yields
that the total spaces
\begin{align*}
&\frac{\Sp(n)}{\SU(p)\times \Sp(l)},~ \frac{\SO(2n+1)}{\SU(p)\times
\SO(2l+1)},~ \frac{\SO(2n)}{\SU(p)\times \SO(2l)},\\&
\frac{\SO(2n)}{\mathbf{S}\tilde\U(n)},~ \frac{\SU(n)}{\SU(p_1)\times
\SU(p_2)}
\end{align*}
of the respective bundles are non-formal under the given conditions.

Due to proposition \ref{propD} (respectively \cite{Oni94}, p.~212), we may replace the stabilisers of
the total spaces of the fibrations by maximal rank subgroups sharing
the maximal torus. We apply this first to the $\SU(p)$-factor which
arose from the fibre-bundle construction and replace it by a
$\mathbf{S}(\U(k_1)\times \dots \times \U(k_t))$. Then we also
substitute the other factors by suitable other ones.

Finally, we are done by an application of proposition
\ref{propA} respectively table \ref{table01}, which allows us
to use the described embeddings of the numerator into larger Lie
groups to form new homogeneous spaces. The result of this process is
depicted in table \ref{HOMtable03}:

The first example we considered here leads to the first line in
table \ref{HOMtable03} by blockwise inclusion $\Sp(n)\hto{}\Sp(N)$ and to
the second line by the inclusion
\begin{align*}
\Sp(n)\hto{}\SU(2n)\hto{}
\SU(N)
\end{align*}
(which is identical to
$\Sp(n)\hto{}\Sp(N')\hto{}\SU(2N')\hto{} \SU(N)$).

The second example produces lines three and four by the inclusions
\linebreak[4]$\SO(2n+1)\hto{}\SO(N)$ and
\begin{align*}
\SO(2n+1)\hto{} \SU(2n+1)\hto{}
\SU(N)
\end{align*}
(which is the same as $\SO(2n+1)\hto{} \SO(N')
\hto{}\SU(2N')\hto{} \SU(N)$).

For examples three and four we do not use any further inclusions.
Example five produces line seven by the inclusion
$\SU(n)\hto{}\SU(N)$.
\end{prf}
The proof guides the way of how to interpret table \ref{HOMtable03}.
Nonetheless, let us---exemplarily for the first line of the
table---shed some more light upon the used inclusions: The inclusion
of the denominator in the first example is given by blockwise
inclusions of
\begin{align*}
\mathbf{S}(\U(k_1)\times \dots \times \U(k_t))\hto{} \Sp(p) \hto{} \Sp(n) \hto{} \Sp(N)\\
\Sp(l_1)\times \dots \times \dots \Sp(l_r) \hto{} \Sp\bigg(\sum_{i=1}^{r} l_i\bigg) \hto{} \Sp(n) \hto{} \Sp(N)\\
\U(l_{r+1})\times \dots \times \dots \U(l_s) \hto{} \Sp\bigg(\sum_{i=r+1}^{s} l_i\bigg) \hto{} \Sp(n) \hto{} \Sp(N)\\
\end{align*}
where $n=p + \sum_{i=1}^s l_i$.

In examples three and four it does not matter whether or not we have
a stabiliser which has---beside the unitary part---special
orthogonal groups of the type $\SO(2l_i)$ only, if there is
additionally one factor of the form $\SO(2l_{r+1}+1)$ or if there is
the factor $\SO(2l_r+1)$ only: All three cases establish maximal
rank subgroups of $\SO(2l+1)$, which therefore share the maximal
torus with $\SO(2l+1)$.

Moreover, note that the relations in table \ref{HOMtable03} are
growing quadratically in $p$. Thus for a constant $l$ we can
always find $p$ and $N$ that will satisfy the restrictions.
\begin{ex}\label{ex02}
Let us describe some of the simplest and most classical consequences
that arise out of our reasoning. We have proved that the spaces
\begin{align*}
\frac{\Sp(n)}{\SU(n)} \qquad\qquad \frac{\SO(2n)}{\SU(n)}
\qquad\qquad \frac{\SU(p+q)}{\SU(p)\times \SU(q)}
\end{align*}
are non-formal for $n\geq 7$, $n\geq 8$ and $p+q\geq 4$
respectively.
\end{ex}
The numerical conditions we impose on the examples seem to be quite sharp. The prerequisites in the last example in table
\ref{HOMtable03} for the space to be non-formal are evidently
satisfied if $p_1,p_2\geq 4$ or if $p_1\geq 3$ and $p_2\geq 5$. These
are sharp bounds:
\begin{prop}\label{prop03}
The homogeneous space
\begin{align*}
\frac{\SU(7)}{\SU(3)\times \SU(4)}
\end{align*}
is formal.
\end{prop}
\begin{prf}
Using the model \eqref{HOMeqn31} from the appendix, we compute a minimal model for the
given homogeneous space as generated by
\begin{align*}
c,~c',~\tilde x &\qquad \deg c=4,~\deg c'=6,~ \deg \tilde x=13
\intertext{with vanishing differentials and by generators representing relations} x_1,~ x_2& \qquad \deg x_1=9,~ \deg x_2=11
\end{align*}
with
\begin{align*}
\dif(x_1)=c\cdot c', ~\dif(x_2)=c^3-(c')^2&
\end{align*}
This model splits as the model of the product of an $F_0$-space and an $\s^{13}$. Thus it is formal.
\end{prf}
As for the interlink between formality and geometry, we remark that,
for example, from the computations in lemma \cite{BK03}.8.2, p.~272,
we directly see that all known odd-dimensional
examples---cf.~proposition\cite{BK03}.8.1, p.~271---of positively
curved Riemannian manifolds are formal spaces; so are the newly announced examples in  \cite{GVZ09} and \cite{PW08}.

The importance of formality issues in geometry is also stressed by
the following result in \cite{BK03}. The non-formality of
\begin{align*}
C:=\frac{\SU(6)}{\SU(3)\times \SU(3)}
\end{align*}
plays a crucial role for the
following (cf.~\cite{BK03}, theorem 1.4): It
gives rise to the existence of a non-negatively curved vector bundle
over  $C\times T$ (with a torus $T$ of $\dim T\geq 2$) which does
not split in a certain sense as a product of bundles over the
factors of the base space, but which has the property that its total
space admits a complete metric of non-negative sectional curvature
with the zero section being a soul.

\begin{rem}
As a generalisation of section \cite{AK11}.4, p.~20, let us recall that it is very easy to construct fibre bundles with interesting properties as far as formality is concerned using non-formal homogeneous space. These bundles illustrate that formality for the base space in general does not have to be related to formality on the total space. Contrast this with a positive result given in theorem \cite{AK11}.A, p.~4.

Homogeneous spaces $G/H$ occur in terms of a fibre bundle
\begin{align*}
H\hto{} G\to G/H
\end{align*}
In the case of a non-formal homogeneous space, this is a fibration
with formal fibre and formal total space, but with non-formal base
space. Moreover, by proposition
\ref{propD} we know that $G/H$ is formal if and only if
$G/T_H$ is formal, where $T_H\In H$ is a maximal torus of $H$. Thus we may successively form the spherical fibre bundles
\begin{align}\label{HOMeqn01}
\s^1\hto{} G/T^k \to G/T^{k+1}
\end{align}
for $k\geq 0$ and $T^k\In T_H$ a $k$-torus. If $G/H$ is non-formal, the space $G/T_H$ is non-formal. Moreover,
we see that there is a choice for $k$ with the property that
\eqref{HOMeqn01} is a fibre bundle with formal total space and
non-formal base space.

Conversely, we may construct a spherical fibre bundle with non-formal total space and formal base space by
\begin{align*}
\s^1 \hto{}  \frac{\SU(6)}{\SU(3)\times\SU(3)} \to \frac{\SU(6)}{\mathbf{S}(\U(3)\times \U(3))}
\end{align*}
The fibration is induced by the canonical inclusion
\begin{align*}
\SU(3)\times\SU(3)\hto{}\mathbf{S}(\U(3)\times\U(3))
\end{align*}
The total space is formal, as it has positive Euler characteristic.
\end{rem}

Let us finally give another example of how to apply theorem \ref{theoA}. This will result in yet another long list of non-formal spaces. In order to shed more light on the motivation for the next theorem---and, in particular, on the inclusions we use---we observe the following: The dimension of the Cartesian product of two spaces clearly is the sum of the dimensions. If both spaces are elliptic, however, the top degree of a non-vanishing rational homotopy group of the product is the maximum over the maximal degrees of non-trivial rational homotopy groups of both factors. So we get the following meta-principle: Taking sufficiently many factors of elliptic Hard-Lefschetz manifolds will result in a Hard-Lefschetz space---with Hard-Lefschetz class the sum of the ones of the factors---to which theorem \ref{theoA} applies. So the base space of the fibration we shall then use will just be the Cartesian product of all the factors, whilst the total space will be gained by ``factoring out'' the Hard-Lefschetz class. This explains, in particular, how the inclusions in the next examples are to be understood.
\begin{theo}\label{theo01}
Every homogeneous space of the form
\begin{align*}
\bigg(\prod_{1\leq i \leq k_1}\Sp(n_{1,i})\times \prod_{1\leq i \leq k_2}\SO(2n_{2,i}+1)\times
\prod_{1\leq i \leq k_3}\SO(2n_{3,i})
\\
\times \prod_{1\leq i \leq k_4}\SO(2n_{4,i})
 \times \prod_{1\leq i \leq k_5} \U(n_{5,i})\bigg)\bigg/ \\
\bigg(\mathbf{S}\bigg(\prod_{1\leq i \leq k_1} \U(p_{1,i}) \times   \prod_{1\leq i \leq k_2}\U(p_{2,i}) \times \prod_{1\leq i \leq k_3} \U(p_{3,i})
\times  \prod_{1\leq i \leq k_4}\tilde\U(n_{4,i})
\\
\times \prod_{1\leq i \leq k_5} \big(\U(p_{5,i})\times \U(q_{5,i})\big)\bigg)\\
\times \prod_{1\leq i \leq k_1}\Sp(l_{1,i})
\times \prod_{1\leq i \leq k_2} \SO(2l_{2,i}+1)\times \prod_{1\leq i \leq k_3}\SO(2l_{3,i})\bigg)
\end{align*}
satisfying
\begin{align*}
&\frac{1}{2}\bigg(\sum_{1\leq i\leq k_1} p_{1,i}^2+4l_{1,i}p_{1,i}+p_{1,i}+
\\&\sum_{1\leq i\leq k_2} p_{2,i}^2+4l_{2,i}p_{2,i}+p_{2,i}+
\\&\sum_{1\leq i\leq k_3} p_{3,i}^2+4l_{3,i}p_{3,i}-p_{3,i}
\\
\geq &  \max\bigg(4 \max_{1\leq i\leq k_1} p_{1,i}+l_{1,i} -1,
\\&4\max_{1\leq i\leq k_2} p_{2,i}+l_{2,i}-1,
\\&4\max_{1\leq i\leq k_2} p_{2,i}+l_{2,i}-5 \bigg)
\end{align*}
is non-formal.
\end{theo}
\vproof

\begin{ex}\label{ex02}
It is obvious that this last method is suited for a variety of generalisations and will produce a myriad of new examples. Let us mention just two more:
\begin{itemize}
\item In remark \ref{HOMrem03} we extend theorem \ref{theoB} to the situation of suitable $\s^3$-fibrations. One source of examples may be a large product of Positive Quaternion K\"ahler manifolds together with the sum of the respective Kraines forms in degree $4$ corresponding to the fibre $\s^3$. As a first example of this kind, we see directly that
    \begin{align*}
    \frac{\Sp(n_1+1)\times \Sp(n_2+1) \times \Sp(n_3+1)\times \dots}{\Sp(n_1)\times  \Sp(n_2) \times \Sp(n_3)\times \dots \times \mathbf{S}(\Sp(1)\times  \Sp(1)\times\Sp(1)\times \dots)}
    \end{align*}
    is a non-formal homogeneous space (provided the number of factors is large). Here, by abuse of notation, $\mathbf{S}(\Sp(1)\times  \Sp(1)\times\Sp(1)\times \dots)$ denotes all the elements $(a_1,a_2,\dots)$ in the quaternionic torus $\Sp(1)\times \Sp(1) \times \dots$  with the property that their product $a_1\cdot a_2 \cdot \dots=1$ in $\hh$. For this we consider $\Sp(1)$ as the unit quaternions.
\item
Using a product involving the Eschenburg biquotient $\s^1\backslash \U(3)/T^2$ it is easy to construct further non-formal genuine biquotients (which are not homogeneous). For this we recall that this biquotient has a Hard-Lefschetz structure, as can be checked easily.
\end{itemize}
\end{ex}
We leave it to the reader to enrich the examples in theorem \ref{theo01} and example \ref{ex02} using theorem \ref{theoA} and proposition \ref{propD}. Doing so one gains even many more non-formal homogeneous spaces/ biquotients.


\section{Non-formal homogeneous manifolds in every large dimension}\label{sec02}

Let us now use the constructed examples of non-formal homogeneous spaces in order to give the
\begin{proof}[\textsc{Proof of theorem \ref{theoE}}]
The example
\begin{align*}
\frac{\SU(N)}{\mathbf{S}(\U(k_1)\times\dots \times \U(k_s))\times
\mathbf{S}(\U(k_{s+1})\times \dots \times \U(k_t))}
\end{align*}
from table \ref{HOMtable03} will serve as a main source of further examples. (It is simply-connected as is the numerator group.)

The restrictions in the table will certainly be
satisfied if $p_1,p_2\geq 4$ or if $p_1\geq 3$ and $p_2\geq 5$ (respectively $N\geq p_1+p_2$). By our computations in appendix \ref{app} we may also use the case $p_1=p_2=3$. We
use a slightly different terminology: Set $p:=p_1$, $k:=p_2$ and
$N\geq 0$. Thus whenever $p=k=3$ or $p,k\geq 4$ or $p\geq 3$ and
$k\geq 5$ the spaces
\begin{align}\label{HOMeqn15}
&\frac{\SU(p+k+N)}{\SU(p)\times
\SU(k)}\\\label{HOMeqn16}&\frac{\SU(p+k+N)}{\SU(p)\times
\mathbf{S}(\U(k-1)\times \U(1))}\\
\label{HOMeqn19} &\frac{\SU(p+k+N)}{\mathbf{S}(\U(p-1)\times
\U(1))\times \SU(k)}\\\label{HOMeqn20}
&\frac{\SU(p+k+N)}{\mathbf{S}(\U(p-1)\times \U(1))\times
\mathbf{S}(\U(k-1)\times \U(1))}\\
\label{HOMeqn25}& \frac{\SU(p+k+N)}{\mathbf{S}(\U(p-2)\times
\U(2))\times\SU(k)}\\
\label{HOMeqn27}& \frac{\SU(p+k+N)}{\mathbf{S}(\U(p-2)\times
\U(2))\times\mathbf{S}(\U(k-1)\times\U(1))}\\
\label{HOMeqn28}& \frac{\SU(p+k+N)}{\mathbf{S}(\U(p-2)\times
\U(2))\times\mathbf{S}(\U(k-2)\times\U(1)\times\U(1))}\\
\label{HOMeqn26}& \frac{\SU(p+k+N)}{\mathbf{S}(\U(p-2)\times
\U(1)\times\U(1))\times\SU(k)}\\
\label{HOMeqn22}& \frac{\SU(p+k+N)}{\mathbf{S}(\U(p-2)\times
\U(1)\times \U(1))\times \mathbf{S}(\U(k-2)\times \U(2))}\\
 \label{HOMeqn21}& \frac{\SU(p+k+N)}{\mathbf{S}(\U(p-2)\times
\U(1)\times \U(1))\times \mathbf{S}(\U(k-2)\times \U(1)\times
\U(1))}
\end{align}
will be non-formal. In the given order the spaces have the following
dimensions:
\begin{align}\label{HOMeqn11}
\nonumber2(k+N)p+(2kN+N^2+1)\\
\nonumber2(k+N)p+(2kN+N^2+2k-1)\\
\nonumber2(k+N+2)p+(2kN+N^2-1)\\
\nonumber2(k+N+2)p+(2kN+N^2+2k-3)\\
2(k+N+4)p+(2kN+N^2-7)\\
\nonumber 2(k+N+4p+(2kN+N^2+2k-9)\\
\nonumber2(k+N+4)p+(2kN+N^2+4k-13)\\
\nonumber2(k+N+4)p+(2kN+N^2-5)\\
\nonumber2(k+N+4)p+(2kN+N^2+4k-13)\\
\nonumber2(k+N+4)p+(2kN+N^2+4k-11)
\end{align}
Regard these homogeneous spaces $M_{p}^{k,N}$ as manifolds parametrised over $\nn$ by the variable $p$.
In other words, we have the infinite sequences \linebreak[4]
$(M^{k,N}_p)^{\eqref{HOMeqn15}}_{p\in \nn}, \dots,
(M^{k,N}_p)^{\eqref{HOMeqn21}}_{p\in \nn}$. We shall now show that
for each congruence class $[m]$ modulo $16$, there is a family
$(M^{k,N}_p)^{m}_{p\in \nn}$ which---for certain numbers
$p$---consists of non-formal manifolds only and which has the
property that
\begin{align*}
\dim (M^{k,N}_p)^{m}_{p\in \nn}\equiv [m] \mod 16
\end{align*}
For this we fix the coefficient of $p$ in each dimension in
\eqref{HOMeqn11} to be $2\cdot 8=16$. The following series realise
the congruence classes $0,\dots,15$. In the third column we give the
smallest dimension for which there is a $p\in \nn$ making the
space non-formal. (From this dimension on the series will produce non-formal examples only.) Once we are given $p$, we may compute this
dimension. We determine $p$ by the following rule: If $k=3$ set
$p=3$, if $k=4$ set $p=4$, if $k\geq 5$ set $p=3$. This guarantees
non-formality. The necessary data is given by the subsequent table:
\begin{center}
\begin{tabular}{ c@{\hspace{10mm}}| @{\hspace{10mm}}c@{\hspace{10mm}}| @{\hspace{10mm}}c}
$[m]$ & series& starting dimension  \\\hline &&\\[-3mm]
$0$& $(M_p^{3,1})_{p\in \nn}^{\eqref{HOMeqn25}}$&  $48$\\[2mm]
$1$& $(M_p^{8,0})_{p\in \nn}^{\eqref{HOMeqn15}}$& $49$\\[2mm]
$2$& $(M_p^{5,1})_{p\in \nn}^{\eqref{HOMeqn26}}$& $50$\\[2mm]
$3$& $(M_p^{4,0})_{p\in \nn}^{\eqref{HOMeqn22}}$& $67$\\[2mm]
$4$& $(M_p^{3,1})_{p\in \nn}^{\eqref{HOMeqn27}}$& $52$\\[2mm]
$5$& $(M_p^{4,0})_{p\in \nn}^{\eqref{HOMeqn21}}$& $69$\\[2mm]
$6$& $(M_p^{3,1})_{p\in \nn}^{\eqref{HOMeqn28}}$& $54$\\[2mm]
$7$& $(M_p^{6,2})_{p\in \nn}^{\eqref{HOMeqn16}}$& $87$\\[2mm]
$8$& $(M_p^{3,1})_{p\in \nn}^{\eqref{HOMeqn21}}$& $56$\\[2mm]
$9$& $(M_p^{4,0})_{p\in \nn}^{\eqref{HOMeqn25}}$& $57$\\[2mm]
$10$& $(M_p^{5,1})_{p\in \nn}^{\eqref{HOMeqn19}}$& $58$\\[2mm]
$11$& $(M_p^{4,0})_{p\in \nn}^{\eqref{HOMeqn26}}$& $59$\\[2mm]
$12$& $(M_p^{7,1})_{p\in \nn}^{\eqref{HOMeqn16}}$& $76$\\[2mm]
$13$& $(M_p^{6,2})_{p\in \nn}^{\eqref{HOMeqn15}}$& $77$\\[2mm]
$14$& $(M_p^{6,2})_{p\in \nn}^{\eqref{HOMeqn20}}$& $78$\\[2mm]
$15$& $(M_p^{6,0})_{p\in \nn}^{\eqref{HOMeqn19}}$& $47$\\
\end{tabular}
\end{center}
\vspace{5mm}
From dimension $72$ onwards every congruence class modulo $16$ can
be realised by one of the given spaces. So we are done.
\end{proof}
Clearly, this theorem is far from being optimal and can easily be
improved by just stepping through table \ref{HOMtable03} or by doing
explicit calculations---see appendix \ref{app}. Indeed, in the proof we have found non-formal
homogeneous spaces in dimensions
\begin{align*}
47,48,49,50,52,54,56,57,58,59, 63,64, 65, 66, 67, 68, 69, 70
\end{align*}
and from dimension $72$ onwards.

We remark that one clearly may use other non-formal series like
\linebreak[4]$\Sp(n)/\SU(n)$ in order to establish similar results. Since the
number of partitions of a natural number $n\in \nn$ is growing
exponentially in $n$, one will find many stabilisers sharing their
maximal tori with $\SU(n)$.

As we pointed out, compact homogeneous spaces are
rationally elliptic.
Thus the question that arises naturally is: In which dimensions are
there examples of non-formal elliptic simply-connected irreducible
compact manifolds?


\section{Proof of proposition D and formality in fibrations}\label{sec03}

Before we shall provide a direct and simple proof of proposition \ref{propD}, let us first see how the result follows from a much more general approach. Using relative obstruction theory we prove in \cite{AK11} (see theorem A, p.~4)
\begin{theo}\label{theo02}
Let
\begin{align*}
F\hto{} E \xto{f} B
\end{align*}
be a fibration of simply-connected topological spaces of finite type. Suppose that   $F$ is elliptic,  formal and   satisfies the Halperin conjecture. Then  $E$ is formal if and only if $B$ is formal.

Moreover, if $B$ and $E$ are formal, then the map $f$ is formal.
\end{theo}
We may apply this to the fibration
\begin{align*}
H/T_H \hto{} G/T_H\to G/H
\end{align*}
which is totally non-cohomologous to zero. Indeed, the rational Leray--Serre spectral sequence of this fibrations degenerates at the $E_2$-term, since homogeneous spaces satisfy the Halperin conjecture. (Recall that this conjecture states exactly this degeneration for arbitrary fibrations of simply-connected spaces with an elliptic fibre of positive Euler-characteristic.)

However, the proof of this result draws on several additional results and can be simplified a lot in our special case. Thus it seems worth while to provide a direct argument.

By $T_H\In H$ we denote the maximal torus. We shall consider the inclusions $i_1\co H\hto{} G$, $i_2\co T_H\hto{} G$ and
$i_3\co T_H\hto{} H$.

\begin{proof}[\textsc{Proof of proposition \ref{propD}}]
The space $G$ admits a minimal Sullivan \linebreak[4]model of the form
$(\Lambda V, 0)$ with $V=V^\textrm{odd}$ of finite dimension; a
minimal model of $H$ is given by $(\Lambda W_H, 0)$ with finite
dimensional $W_H=W_H^\textrm{odd}$. Equally, we have a minimal model
for $T_H$ given by $(\Lambda W_T, 0)$ with finite dimensional
$W_T=W_T^1$. Thus a model of $G/H$ is given by $(\Lambda
W_H^{+1}\otimes \Lambda V,\dif_1)$ and a model for $G/T$ is
given by $(\Lambda W_T^{+1}\otimes \Lambda V,\dif_2)$---see \cite{FHT01} p.~218. (By the superscript $+1$ we denote a corresponding degree shift, as usual.)

Both in the case of $G/H$ as well as for $G/T_H$ we see that the
corresponding Sullivan models are pure algebras: The differential
$\dif\in \{\dif_1,\dif_2\}$ in each case has the property that $\dif|_{\Lambda
W^{+1}}=0$ (as $(\Lambda W^{+1},0)$ forms a differential
subalgebra) and that $\dif|_{\Lambda V}\in \Lambda W^{+1}$ by
\cite{FHT01} p.~217. This directly yields a pure filtration. Clearly, the
corresponding minimal models then are pure, too.

The differentials are given by $\dif_1 v=H(\B i_1)(v^{+1})$ (for $v\in V$ and in the case of $G/H$) respectively by $\dif_2 v=H(\B i_2)(v^{+1})$ (for $G/T_H$). Since $H^*(\B H)=H^*(\B T_H)^{\W(H)}$, i.e.~the cohomology of the classifying space of a Lie group consists of the invariant polynomials in the cohomology of the classifying space of the maximal torus under the action of the Weyl group, we derive the following: The map $H^*(\B i_3)$ induced by the inclusion $i_3\co T_H\hto{} H$ is also just the inclusion in this setting. Consequently, on the level of differentials we obtain
\begin{align}\label{eqn07}
\dif_2 v=H(\B i_2)(v^{+1})=H(\B (i_3\circ i_1))(v^{+1})=H(\B i_3)(\dif_1 v^{+1})
\end{align}
Thus we may identify the differentials $\dif_1$ and $\dif_2$ from the models for $G/H$ and for $G/T_H$. This will allow us to deduce that one space is formal if and only if so is the other:

\vspace{5mm}

First assume $G/T_H$ to be formal. By lemma \ref{lemma02} and theorem \cite{FHT01}.14.9 we obtain the decomposition
\begin{align}\label{eqn05}
(\Lambda W_{T_H}^{+1}\otimes V,\dif_2) \cong (\Lambda V',\dif) \otimes (\Lambda \langle z_1, \dots , z_l\rangle,0) \otimes (C,\dif)
\end{align}
where $(\Lambda V',\dif)$ is minimal of positive Euler characteristic and where $(C,\dif)$ is a contractible algebra. Due to the formality of $G/T_H$ we obtain that $(V')^\even=(V')^\odd$. Since we may identify the differentials for $G/H$ and $G/T_H$ as above, we derive a similar splitting
\begin{align}\label{eqn06}
(\Lambda W_H^{+1}\otimes V,\dif_1) \cong (\Lambda V'',\dif) \otimes (\Lambda \langle z'_1, \dots , z'_{l'}\rangle,0) \otimes (C',\dif)
\end{align}
Indeed, by the constructions in the proofs of lemma \ref{lemma02} and theorem \cite{FHT01}.14.9, the right hand side is contained in the left hand side and the isomorphism is rather an equality. Thus decomposition \eqref{eqn06} can be established in analogy to decomposition \eqref{eqn05} considering the following arguments: It might happen that the image $\dif_1(a)$ of an element $a\in V$ has an element $b$ from $W_H^{+1}$ as a non-trivial summand. In this case, both elements $a,b$---with $\deg a$ being odd and $\deg b$ being even---may be taken to lie in $C'$. (In particular, the element $a$ does not lie in $V''$.) Thus the dimension of $V''$ may be smaller than the dimension of $V'$. However, the Euler characteristic is preserved, i.e.~$\chi(V'')=\chi(V')=0$.

Note that the free part $(\Lambda \langle z_1, \dots , z_{l}\rangle,0)=(\Lambda \langle z'_1, \dots , z'_{l'}\rangle,0)$ is not altered, since an element with non-trivial differential in $(\Lambda V',\dif)$ has non-trivial differential in \eqref{eqn06} and vice versa due to \eqref{eqn07}.

Consequently, we have that $\deg (V'')^\even=\deg (V'')^\odd$ and the minimal model of $G/H$ splits as in lemma \ref{lemma02}. In other words, the space $G/H$ is formal.

Conversely, given the formality of $G/H$ and the splitting \eqref{eqn06}, the decomposition \eqref{eqn05} can be derived by an analogous inverse process again preserving the Euler characteristic $\chi(V'')=\chi(V')=0$. This yields the formality of $G/T_H$.
\end{proof}

\begin{cor}\label{cor01}
A biquotient $\biq{G}{H}$ is formal if and only if $\biq{G}{T_H}$ is
formal.
\end{cor}
\vproof

\vspace{5mm}

There are several questions which may be stated to further explore the Halperin conjecture in the context of formality and ellipticity. One is
\begin{ques}[Lupton \cite{Lup98}, p.~18]\label{ques01}
Let $F\hto{} E\to B$ be a fibration in which $F$ is formal and
elliptic and $B$ is formal. If $E$ is formal, then is the fibration
totally non-cohomologous to zero?
\end{ques}
A partial positive answer in the case of odd dimensional spheres as a base is given in
proposition \cite{Lup98}.3.3, p.~11.

Although, in general, this question can be trivially answered in the negative via the example of the universal fibration $G\hto{} \E G \to \B G$ where $G$ is a non-trivial Lie group together with its classifying space $\B G$,
it is certainly worth while to impose further conditions on $B$, for example, and to investigate this question and related ones in this context. However, the question can also be answered in the negative if ellipticity on the base space is imposed. Here the Hopf bundles $\s^1\hto{} \s^3 \to \s^2$ or $\s^3\hto{} \s^7 \to \s^4$ (for the simply-connected case) may serve as an example.

Any positive answer to question \ref{ques01} might be considered as a partial converse to theorem \ref{theo02} in a certain sense.


\appendix
\section{A computational example}\label{app}

Let us illustrate a
showcase computation by which we reprove the non-formality of
$\frac{\SU(6)}{\SU(3)\times\SU(3)}$ via a direct computation.
This serves the following purposes:
\begin{itemize}
\item We illustrate how direct computations as in the proof of proposition \ref{prop03} work.
\item We show that any sort of sufficiently ``small'' respectively low-dim\-ension\-al example can easily be provided. And most important of all:
\item We illustrate how fast the complexity of computations increases making it absolutely necessary to provide general principles as we did.
\end{itemize}

We construct a model for the homogeneous
space according to theorem \cite{FOT08}.3.50. A minimal Sullivan
model for the Lie group $\SU(n)$ is given by $\Lambda (\langle x_2, x_3, \dots, x_n\rangle,0)$ with $\deg
x_i=2i-1$; a Sullivan model for the classifying space
$\B\SU(n_k)$ (with $n_k\geq 2$) is given by the polynomial
algebra $(\Lambda \langle c_2,\dots,c_{n_k}\rangle,0)$
with $\deg c_i=2i$. We obtain a model
\begin{align}\label{HOMeqn31}
\nonumber&\APL\bigg(\frac{\SU(n)}{\SU(n_1)\times\dots\times\SU(n_k)},\dif\bigg)\\\simeq
&\Lambda(\langle c_2^1,\dots, c_{n_1}^1, c_2^2,\dots,
c_{n_2}^2,\dots, c_2^k,\dots, c_{n_k}^k\rangle\otimes \langle
x_2,\dots x_{n}\rangle,\dif)
\end{align}
with $\deg c_i^{j}=2i, \deg x_i=2i-1$. The differential vanishes on
the $c_i^j$ and \linebreak[4]$\dif x_i=H^*(\B(\phi))y_i$, where $y_i$ is the
generator of $\B\SU(n)$ corresponding to $x_i$ and $\phi$ is the
blockwise inclusion map
$\SU(n_1)\times\dots\times\SU(n_k)\hto{}\SU(n)$.

Note that we may take the $c_i^j$ for the $i$-th universal Chern
classes of $\B\SU(n_j)$. Thus blockwise inclusion of the $\SU(n_j)$
yields the following equations
\begin{align*}
\dif(x_2)&=\sum_{i=1}^k c_2^i\\
\dif(x_3)&=\sum_{i=1}^k c_3^i\\
\dif(x_4)&=\sum_{i=1}^k c_4^i + \sum_{i\neq j} c_2^ic_2^j\\
&\vdots\\
\dif(x_i)&=(c^1\cdot \dots \cdot c^k)|_{2i} \\
&\vdots\\
\dif(x_n)&=c_{n_1}^1\cdot \dots \cdot c_{n_k}^k
\end{align*}
which become obvious when writing the Chern classes as elementary
symmetric polynomials in elements that generate the cohomology of
the classifying spaces of the maximal tori. By $c^j$ we denote the
total Chern class $1+c_2^j+c_3^j+\dots + c_{n_j}^j$ of $\SU(n_j)$
and by $|_{2i}$ the projection to degree $2i$.

In the case $k=2$, $n_1=n_2=3$ we obtain a model $(\Lambda V, \dif)$, where the graded vector space $V$ is generated by
\begin{align*}
\begin{array}{ll}
c,c',d,d'  &\deg c=\deg d=4,~ \deg c'=\deg d'=6 \\
x_2,x_3,x_4,x_5,x_6 &\deg x_2=3,~ \deg x_3=5,~ \deg x_4=7,~ \deg x_5=9,~\\&
\deg x_6=11
\end{array}
\end{align*}
The differential is given by
\begin{align*}
\dif c&=\dif c'=\dif d=\dif d'=0\\
\dif x_2&=c+d,~ \dif x_3=c'+d',~ \dif x_4=c\cdot d,~ \dif x_5=c'\cdot
d+d'\cdot c,~\dif x_6=c'\cdot d'
\end{align*}
Thus the cohomology algebra of this algebra is given by generators
\begin{align*}
e,e',y,y' &\qquad\deg e=4,~ \deg e'=6,~\deg y=13, ~\deg y'=15
\end{align*}
and by relations
\begin{align*}
e^2=e\cdot e'=e'^2=y\cdot e=y^2=y\cdot y'=y'^2=y'\cdot e'=0,\quad
y\cdot e'=-y'\cdot e
\end{align*}
Hence Betti numbers are as follows:
\begin{align*}
b_0=b_4=b_6=b_{13}=b_{15}=b_{19}=1
\end{align*}
with the remaining ones equal to zero.

Using the construction principle from \cite{FHT01}.12, p.~144--145, we compute the
minimal model for this algebra. On generators it is given by
\begin{align*}
e,e',x,x',x'' \qquad \deg e=4,~\deg e'=6, ~\deg x=7, ~\deg x'=9,~ \deg
x''=11
\end{align*}
The differential $\dif$ is given by
\begin{align*}
\dif e=\dif e'=0, ~\dif x=e^2, ~\dif x'=e\cdot e', ~\dif x''=e'^2
\end{align*}
It is easy to check that the Massey product $\langle e,e,e'\rangle$ is non-trivial. So $\frac{\SU(6)}{\SU(3)\times\SU(3)}$ is not formal---see proposition \cite{FOT08}.2.90, p.~94.

We remark that the spaces $\SU(4)/(\SU(2)\times \SU(2))$ and
$\SU(5)/(\SU(2)\times \SU(3))$ are formal as a straight-forward
computation shows.




\vfill

\begin{center}
\noindent
\begin{minipage}{\linewidth}
\small \noindent \textsc
{Manuel Amann} \\
\textsc{Department of Mathematics}\\
\textsc{University of Pennsylvania}\\
\textsc{David Rittenhouse Laboratory}\\
\textsc{209 South 33rd Street}\\
\textsc{Philadelphia, Pennsylvania}\\
\textsc{PA 19104-6395} \\
\textsc{USA}\\
[1ex]
\textsf{mamann@uni-muenster.de}\\
\textsf{http://hans.math.upenn.edu/$\sim$mamann/}
\end{minipage}
\end{center}

\end{document}